     \documentclass{article}
     \usepackage[english]{babel}
     \usepackage[leqno]{amstex}
     \usepackage{theorem}
     \usepackage[bfhd]{lifodi} % Private Package
     \usepackage{varioref,xspace}
     \usepackage{program}
\begin{document}
%% \BEGIN{lifo00} %% Macros, Initialisierungen, Abstract.
%
%   LIFO00.TEX  %  Vorspann nach \begin{document}
%--------------------  Makros
    \let\Include\include
    \let\include\input
    \makeatother
    \hfuzz=22pt %% overfull box < 22pt wird nicht gemeldet
%================================================================
\title{% \vspace*{22mm} %{20mm}%{-10mm}
    \Large\bf
    A LiFo Dynamic Dictionary
    }
\author{\large Sch\"onbrunner Josef %% \quad \today \quad Version v1.0
    \\
    \normalsize
    Institut f\"ur Logistik der Universit\"at Wien
    \\
    \normalsize
    Universit\"atsstraae 10/11,
    A-1090 Wien (Austria)
    \\
    \normalsize
    e-mail a8121dab@@helios.edvz.univie.ac.at
} \date{}
\maketitle

\begin{abstract}
    % LIFO.ABS
    % ~~~~~~~~
    Data structures that realize a dictionary are characterized
    by three \iftrue basic instructions \else fundamental operations\fi:
    (1) \textit{Insert} a new entry ({\em key} and {\em value}).
    (2) \textit{Search} by a key, returning the associated value.
    (3) \textit{Delete} an entry.
    %\end{itemize}
    \par
    Known realizations are hashing schemes and various types of search trees.
    Time complexity of the fundamental operations is measured as a function of
    the number of entries $n$, the (binary) size of a key $s$
    is usually not considered.
    For search trees the expected time as well as the upper limit of time
    is $O(\log n)$.
    \textit{LiFo dictionary} is a new implementation, the
    time limits for \textit{Insert} and for \textit{Search} both are a linear
    function of the length $s$ of the used key, that is $O(s)$.
    The \textit{LiFo dictionary} furthermore provides two additional
    basic operations: (4) \textit{open environment}
    and (5) \textit{close environment} with a constant time for performing.
    \textit{Close environment} needs only one assignment by which it
    restores exactly the same internal situation
    as before the last call of \textit{open environment}.
    This feature cannot be realized by any of the data structures
    mentioned above, therefore the prefix LiFo (=\textit{last in first out}),
    another name for \textsf{stack}.
    This ability is highly suitable for software applications that
    frequently perform local symbol to value binding throughout multiple
    levels of environments.
    E.g. \textit{Lisp-Interpreter} or  D. Knuth's \TeX programm are such
    applications.
\end{abstract}

\vspace{1\baselineskip} %% {1ex}
\noindent
\textbf{Keywords: }  Datastructure Dictionary, LiFo-Facility,
              Computational Comlexity.
% -----------------------------------------------+
%% \END{lifo00}

%% \BEGIN{lifo01} % LIFO01.TEX
%                   ~~~~~~~~~~
\section{Introduction}

To illustrate the principle of this data structure in a simplified
fashion, a game model, which is demonstrated by an example,
seems to be suitable.
Assume that the president of some association keeps a dictionary of
the members.
Keys of the dictionary are the member's names, information to be accessed
as value assigned to a key is the year of birth.
Assume the current state of the dictionary are three entries:
ALFRED (1940), ALBERT (1955), PETRA (1960).
This current state in the game is modelled by a route pattern
of square shaped cards, each with one letter written on it,
as the following figure shows.

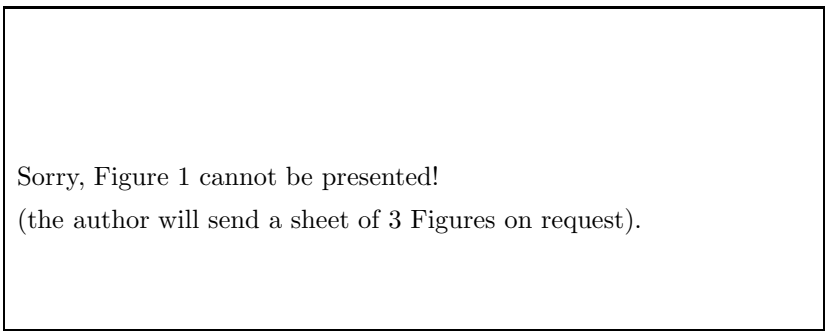
\Nofigbox{.9}{-0.6in}{1.6in}{%
    \begin{minipage}{.87\textwidth}
        Sorry, Figure 1 cannot be presented!
        \\[1ex]
        (the author will send a sheet of 3 Figures on request).
    \end{minipage}
}
%\vspace*{1.7in}
%\IfFileExists{fig1.bmp}{?}{\nofigbox{.9}{-0.6in}{1.6in}{1}}

We label card \framebox[1.5em]{A} by $\triangleright$ as starting square
and mark the ends with the value  encircled (the birth year). % elliptically.
Using another set of these cards beside we place a so called argument line,
whose letters form the word to be searched (the input).
The search is simulated by the following procedure:
First we put a piece \textit{a} onto the starting square
of the above route pattern and a piece \textit{b} onto
the first square of the argument line. This is the starting configuration.
A configuration of the game in process is characterized by the square
locations of the two pieces, we shall refer to it as the current squares
in the dictionary and in the argument line.
A successful final configuration is reached, if \textit{a} covers the
last square of the line, \textit{b} is upon a square with circled
value attached and the letters of both squares agree.
The answer of the search then is the circled value next to \textit{b}.
The change from one configuration to the next is described as follows:

\begin{quote}
    If the letters of the current fields are the same, the pieces
    move one square in straight direction.

    In the other case, if on the left hand of piece \textit{b}
    a square is adjacent, then \textit{b} moves left to this square
    and \textit{a} keeps its position; if there is no square lefthand,
    however, the game fails, that is the argument does not appear in
    the dictionary.

    If instead of \textit{searching} an extension of the dictionary
    is to be performed by \textit{inserting} the input as a key,
    then we continue in the last case as follows.
    We place a card with the same letter as that below piece \textit{a}
    so that it adjoins lefthand to the sqare where piece \textit{b}
    resides. We twist the card a small angle to the left, then we append
    a copy of the part of the argument line in front of piece \textit{b}
    as a new branch. To finish \textit{inserting}, a further argument,
    carrying the value to be associated with the key is required.
    (The value is the encircled year in the example above)
\end{quote}

If we extend our example by inserting an entry for PETER (1965) we obtain:

\nofigbox{.9}{-0.6in}{1.6in}{2}
%%\vspace*{1.7in}

This game should illustrate only how the mechanism of the
data structure works in principle.
We observe the following drawback in this model:
From a certain size upwards considerable multiple branching will arise
so that search becomes inefficient. A simple rule helps to overcome
this difficulty: rely on a binary coding of the letters and play the game
with the alphabet of 0 and 1!
Let us record that access by a matching chain of bits (input strip)
is realized by a sequence of moves in the game model.
Each move comprises \textit{testing agreement} of the two current
(the focused)
bits of the input strip and the search tree respectively;
if testing results to \textit{false}, the focused bit migrates
to the leftside branch (which exists for a matching input)
and \textit{testing agreement} is repeated; if finally testing results
in \textit{true}, both focuses move one node ahead (to the next bit).
For a correct input the focus migrates through the branch in the
search tree that carries the same sequence of bits as the input strip.
It is significant to record that each entry except the first (the oldest one)
is mapped to exactly one branching node in the search tree.
Now let us design an adequate data structure.
%% \END{lifo01}

%% \BEGIN{lifo02}
% LIFO02.TEX
% ~~~~~~~~~~
\section{Data Structure of \LEXIKON}

An entry \textit{key + value} is recorded as a data type
\lexkonf with field components. % (alias field selectors)
The dictionary as a whole is an organized collection of objects
of type \lexkonf together with a pool for keys, that is a compound
of objects, which are interconnected by references, sometimes called
pointers.
One component to be considered is the bit-string of the key of an entry.
With exception of the first entry the right part, that is
the continuation after a branching off from the key of an elder entry,
is sufficient for the searching algorithm.
We shall record this bit string as a component of \lexkonf
referred to by selector name \rkey.
Now let us imagine further branching off nodes located in certain
places of \rkey. Remember the fact that each entry
is mapped to exactly one branching node in the search tree, namely
the one which leads into the entry. Hence it seems to be advisable
to record only the first branching off by a component attached to
the structure of \lexkonf and to delegate further branching offs
to the entry reachable from the first one.
This is managed by introducing three further components for \lexkonf:
\textbf{\dist} as
the place of first branching in bit string \rkey;
\textbf{\zweig} as pointer to the entry (object of type \LEXKONF)
that is reachable by the first branching off;
and finally
\textbf{\link}
as pointer to record further branching offs,
which belong to the bit string \rkey of the entry that was the \fokus
(the current \LEXKONF)
before the last assignement $\fokus \gets \fokus\REF\bdot\zweig$.

A pointer or reference to a certain object is the address or position
of the object within the digital storage.
A declaration \REF\lexkonf \textit{lexpos} stipulates \textit{lexpos}
as a pointer to an object of type \lexkonf, the object is denoted
by \textit{lexpos}\REF.

The meaning of the components \zweig and \link is illustrated
by the following sketch that enters into our game model.
Assume the dictionary contains 4 entries, the adresses of these
are denoted by $L_1,L_2,L_3,L_4$ and the corresponding route pattern is:

%\vspace*{1in}
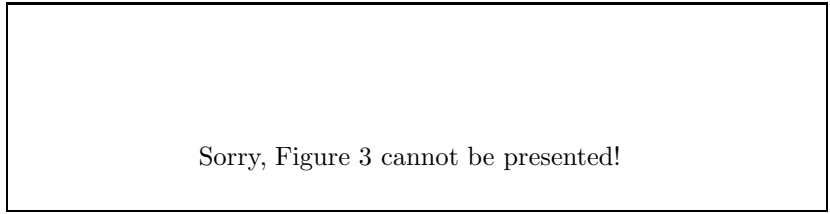
\nofigbox{.9}{-0.2in}{1in}{3}

Then we have
\begin{math}
    L_2=L_1\REF\bdot\zweig, \quad
    L_3=(L_1\REF\bdot\zweig)\REF\bdot\link,
    \quad \text{and} \\
    L_4=((L_1\REF\bdot\zweig)\REF\bdot\link)\REF\bdot\link
\end{math}.
If for an object $X$ of type \lexkonf there are $n+1$ entries
corresponding to the branching offs on bit string $X\bdot\rkey$
according to the above model,
that is each key associated with such an entry
(component \bdot\rkey) exhibit the first difference
to $X\bdot\rkey$ in the branching locations,
then the entries are
\begin{math}
    X\bdot\zweig\REF, X\bdot\zweig\REF\bdot\link\REF,
    \dots,
    X\bdot\zweig\REF(\bdot\link\REF)^n
\end{math}.
The branching positions can be read out of these objects:
component $Y\bdot\dist$ for each
$Y=X\bdot\zweig\REF(\bdot\link\REF)^j\quad(j\le n)$
is the information that the leftmost difference
of the key of $X$ and that of $Y$ is the bit at position $Y\bdot\dist$
in $X\bdot\rkey$
(although \rkey in general is only a right piece of the real key).

\begingroup %>>>>>>>>>>>>>>>>
%    \selectlanguage{english}

\begin{thms}{Preliminary Data Types}
    \begin{program}
        \INDENT %
        \deftype \BIT = \encurs{0,1};\qquad %
        \deftype \BOOL = \encurs{\false,\true};
        \deftype \BITSTR = \ARRAY[\star]\quad \OF \BIT;
        \comment %
            \BITSTR s \quad\textbf{means} \quad%
            s = \enangle{s[i]}_{i<\length(s)}%
        \endcomment
~
        \deftype \NATN=\text{{Elements of $\mathbb{N}$}};
        \deftype \VALUE=\text{{any type you may need below}};
        \deftype \LOC=%
        \text{position of storage referring to some object};
        \UNDENT
    \end{program}
\end{thms}

\begin{Thms}{Data Types and References}\em
    For each data type \textsc{T}\ we consider the type \REF\textsc{T}\
    of reference to an object of type \textsc{T}\
    (location of the object within storage).
    If a variable $p$ is declared to be of type \textsc{T}\
    (usually by writing \inanf{$\textsc{T}\ p$}) then $p\REF$ denotes
    the object referred to by $p$.
    The relationship between \REF\textsc{T}\ and \LOC\ is that
    under declarations
    \inanf{$ \textsc{T}\ p,q; \LOC\ r; $}
    the commands
    \inanf{$ r \gets p; q \gets r; $} yield $p=q$.
\end{Thms}

\begin{thms}{Procedures of Storage-Allocation}\em
    \begin{program}
        \INDENT %
        \REF\textsc{T}\ \valued %
        \mbox{/* \textsc{T} may be any type */}
        \create(\textsc{T}, \var \LOC\ \textit{freeloc}) \idgl
        \mbox{/* \textit{freeloc}\ is considered as the top of a stack */}
        \BEGIN
            \REF\textsc{T}\ \textit{result} \gets \textit{freeloc};
            \text{Increase \textit{freeloc}\ by the size of \textsc{T}};
            \RET(\textit{result});
        \END.
~
        \REF\textsc{T}\ \valued %
        \mbox{/* \textsc{T} may be any type */}

        \ccopy(\textsc{T}\ v , \var \LOC\ \textit{freeloc}) %
        \idgl
        \BEGIN
            \REF\textsc{T}\ \textit{result} \gets \textit{freeloc};
            \text{Increase \textit{freeloc}\ by the size of \textsc{T}};
            \textit{result}\REF \gets v;
            \RET(\textit{result});
        \END.
    \end{program}
\end{thms}

\par
First we study a simpler version of the data structure \lexkonf,
which does not yet realize the LiFo-property but differs from the
final version essentially only with regard to the operation
of adding a new entry by procedure \textit{insert}.
The logical connectives $\land$ and $\lor$ are used with
semantics that differ from the one used in logic, as
they will be interpreted as a \textit{sequential and} and \textit{or}
respectively. That is, if the first operand of $\lor$ yields
\true\ within an evaluation, then evaluation of the expression
built upon $\lor$ ceases with value \true\ without any computation
of the second operand. In a similar way this applies to
$\land$ with \false\ in place of \true.
This can prevent e.g. for a conditional instruction as
$\IF (\findpos \neq \nil \land d \gt \findpos\REF\bdot\dist) \dots$
the computation of an undefined expression $\findpos\REF\bdot\dist$
in case of $\findpos=\nil$.

\NumberProgramstrue
\begin{DATASTR} \DSLEX \quad \mbox{/* simple (non-lifo) Version */}
    \label{lexikon}
    \begin{program}
        \INDENT %
        \deftype \lexkonf = \struct
        \BEGIN
            \REF\BITSTR \rkey;
            \NATN \dist;
            \REF\lexkonf \zweig,\link;
            \VALUE\yspace \val;
        \END
~
        \vspace{1\baselineskip}
        \deftype \LEXIKON=
        \struct \BEGIN \REF\lexkonf \gate; \LOC\ \freeloc; \END;
~
    \begin{minipage}^{\noobeycr}[t]{.93\textwidth}
        /*
        Procedure \procname{search} computes the
        \REF\LEXIKON position $\findpos$ for a bitstring $x$
        that corresponds to an entry in the dictionary $L$,
        the key of which is $x$, provided that this key
        occurs in $L$.
        The proper answer after a successful search
        then is \findpos\REF\bdot\val.
        The answer whether  the search was successful
        is the return value of procedure \procname{search}.
        Procedure \procname{insert} relies on \procname{search}.
        If $\procname{search}(\dots)$ returns value \true,
        then \findpos\REF\bdot\val needs only to be rewritten
        by the new associated value (a parameter of \procname{insert}).
        In the other case $L$ will be enlarged by one element
        (type \lexkonf). Particularly for this purpose \procname{search}
        is provided with output parameters \insertpos\ and $d$.
        */
    \end{minipage}
~
    \BOOL \enspace \valued
    \procname{search}%
    (\LEXIKON L,\ \BITSTR x,
      \hspace{12mm}%
      \var \REF\lexkonf \findpos,\insertpos,\ \var \NATN $d$) \idgl
    \BEGIN
        \findpos \gets L\bdot\gate;
        \label{r1}%
        \IF (\findpos = \nil) \THEN \RET(\false); %
        \mbox{/* empty dictionary */} \untab \untab
        \WHILE(\true) \DO\quad \mbox{/* terminated by \RET */} \untab
        \BEGIN
            \IF (\findpos\REF\bdot\rkey\REF = x) %
            \THEN \RET(\true); \untab\untab
            \IF %
             (\text{ one of } x \text{ and } \findpos\REF\bdot\rkey\REF %
                \text{ is a left substring of the other }) \untab
            \THEN d \gets %
                \min(\length(x),\length(\findpos\REF\bdot\rkey\REF)) + 1;
            \ELSE \keyword{let}\yspace d \text{ be the smallest } d, %
                \text{ such that }%
                x[d] \neq (\findpos\REF\bdot\rkey\REF)[d]; \untab
~
            \insertpos \gets \findpos; %
            \mbox{/* $\neq\nil$ %
                (on account of \textit{\ref{r1}}\ and \textit{\ref{r2}}) */}
            \findpos \gets \findpos\REF\bdot\zweig;
            \WHILE (\findpos \neq \nil \land d \gt \findpos\REF\bdot\dist) %
            \DO \untab
            \BEGIN
                \insertpos \gets \findpos;
                \findpos \gets \findpos\REF\bdot\link;
            \END;
            \label{r2}% \Delta:\enspace %
            \IF (\findpos = \nil \lor d \neq \findpos\REF\bdot\dist)
                \THEN \RET(\false); %
            \mbox{/* $x$ is not a valid key */} \untab
            \UNDENT
            x \gets \enangle{x[d+j]}_{j<\length(x)-d}; %
            \mbox{/* $\idgl \textsf{shift_right}(x,d)$ */}
        \begin{minipage}^{\noobeycr}[t]{.84\textwidth}
            /* search continues at the entry the branch off
                leads to, i.e. after position $d$ of $\rkey\REF$.
             */
        \end{minipage}
        \END;
    \END. \mbox{/* \textit{search} */}
~
    \procname{insert}(\LEXIKON L,\ \BITSTR x,\ \VALUE\ v) \idgl
    \BEGIN
        \LEXKONF\ \findpos, \insertpos; \quad \NATN d;
        \IF (\procname{search}(L,x,\findpos,\insertpos,d)) \untab
        \THEN \findpos\REF\bdot\val \gets v
        \ELSE \untab
        \BEGIN
            \newpos \gets \create(\LEXKONF, L\bdot\freeloc);
            \newpos\REFD\dist \gets d;
            \IF (\insertpos = \nil) \THEN L\bdot\gate \gets \newpos; %
            \untab \untab
            \mbox{/* note that } \insertpos=\nil %
            \mbox{ iff } L\bdot\gate = \nil \mbox{ */}
            \IF (\insertpos\REFD\branch = \findpos) \untab
            \THEN \insertpos\REFD\branch \gets \newpos
            \ELSE \insertpos\REFD\link \gets \newpos; \untab
            \newpos\REFD\branch \gets \nil;
            \newpos\REFD\link \gets \findpos;
            \newpos\REFD\rkey \gets %
            \ccopy(\enangle{x[d+j]}_{j<\length(x)-d}, L\bdot\freeloc);
            \newpos\REFD\val \gets v;
        \END; \untab
    \END. \mbox{/* \textit{insert} */}
\end{program}
\end{DATASTR}

\begin{observ}\em
    The worst case time complexity of procedure \procname{search}
    of \ref{lexikon} is $\mathsf{O}(s)$ where $s=\textit{length of } x$.
\end{observ}

\begin{proof}
    The instructions within the while loop at lines \tit{21-40}
    are to be considered as follows: Assume the loop is passed
    $r$ times, the associated values of $d$ being $\enangle{d_j}_{j<r}$.
    Looking at lines \tit{25,26} and \tit{38} we have
    \begin{math}
        d_0 + d_1 + \dots + d_{r-1} \le s=\length(x)
    \end{math}.
    More than two subsequent $d_i$'s cannot be  equal to 0
    ($d_i=d_{i+1}=0$ is only if keys $u, u0v, u1w$ occur),
    therefore $r \le 2\dot s$.
    As to the inner loop of lines \tit{30-34},
    the construction of the \link\ components
    according to \procname{insert} ensures that
    $\tit{anypos}(\REFD\link)^n\REFD\dist$ is a strictly increasing
    sequence limited by $d$.
    In this way the inner loop is repeated $p_j$ times, if
    \begin{displaymath}
        0 \le \findpos\REFD\dist < \findpos\REFD\link\REFD\dist
        < \dots <
        \findpos(\REFD\link)^{p_j - 1}\REFD\dist
        \le d_j
    \end{displaymath}
    If we take
    \begin{math}
        d_{j\ k}=
        \findpos(\REFD\link)^{k+1}\REFD\dist
        - \findpos(\REFD\link)^{k}\REFD\dist
    \end{math}
    then each $d_{j\ k} \ge 0$ and
    $d_j = d_{j\ 0} + \dots + d_{j\ p_j - 1}$.
    If we concatenate the tuples $\enangle{d_{j\ k}}_{k<p_j}$
    and let $\enangle{\delta_i}_{i<h}$ be the result of the concatenation,
    we conclude, as above for the sequence of the $d_j$'s,
    that $h \le s$, where $h$ is the sum of the number of repetitions
    of the small loop for each traverse of the big loop.
    Hence the amount of computation time due to the small loop
    is $\mathsf{O}(s)$, the number of repetitions of the big loop
    is $\mathsf{O}(s)$, the instructions outside the small loop
    either require constant time or are limited by $d_j$ and, as
    the sum of the $d_j$'s is not greater then $s$, the whole algorithm
    of \procname{search}\ has worst case time complexity $\mathsf{O}(s)$.
\end{proof}
\endgroup %%%<<<<<<<<<<<<<<<<
%% \END{lifo02}

%% \BEGIN{lifo03}
% LIFO03.TEX
% ~~~~~~~~~~
\NumberProgramstrue
\begin{DATASTR}\em \DSLLEX \label{lifolex}
\\ \hspace*{10mm}
    \begin{minipage}[t]{.84\textwidth}
        % Inherit all of \LEXIKON except \procname{insert}
        Each component of the preceding data-structure is included,
        partly modified.
        Two additional procedures \procname{open_environment}\
        and \procname{close_environment}\ are introduced.
        Data-type \LEXIKON\ gets the new component \fieldsel{pop}
        of type \REF\LEXIKON, so that a stack of \LEXIKONnsp-frames
        is maintained.
    \end{minipage}
\begin{program}
    \INDENT %
    \deftype \lexkonf = \mbox{as in \ref{lexikon}\ \LEXIKON};
~
    \deftype \LEXIKON=
    \struct \BEGIN
        \REF\lexkonf \gate;
        \LOC\ \freeloc;
        \REF\LEXIKON \fieldsel{pop}
    \END;
~
    \procname{open_environment}(\var \REF\LEXIKON \dictptr) \idgl
    \BEGIN
        \REF\LEXIKON \textit{Old}\dictptr \gets \dictptr;
        \LOC\ \textit{NewFreeLoc} \gets \dictptr\REFD\freeloc;
        \dictptr \gets \create(\LEXIKON, \textit{NewFreeLoc});
        \dictptr\REFD\gate \gets %
        \ccopy(\textit{Old}\dictptr\REFD\gate\REF, \textit{NewFreeLoc});
        \dictptr\REFD\fieldsel{pop} \gets \textit{Old}\dictptr;
        \dictptr\REFD\freeloc \gets \textit{NewFreeLoc};
    \END.
~
    \procname{close_environment}(\var \REF\LEXIKON \dictptr) \idgl
    \qquad    \dictptr \gets \dictptr\REFD\fieldsel{pop}
~
    \begin{minipage}^{\noobeycr}[t]{.87\textwidth}
        /*
        Note, that we associate our dictionary with a variable
        \REF\LEXIKON \dictptr. The \textsf{LiFo}-storage for the whole
        compound is characterized by its top-position \dictptr\REFD\freeloc.
        Thus the single assignment of \procname{close_environment}
        suffices to restore the state before the last call of
        \procname{open_environment}.
        */
    \end{minipage}
~
    \BOOL \enspace \valued
    \procname{search}(\REF\LEXIKON \dictptr,\ \BITSTR x) \idgl
    \BEGIN
        \REF\lexkonf \findpos,\insertpos; \NATN $d$;
        \mbox{/* these are no longer required for \procname{insert}, */}
        \mbox{/* we use \procname{search} of \ref{lexikon} */}
        \RET(\DSLEX\bdot\procname{search}%
            (\dictptr\REF, x, \findpos,\insertpos));
    \END.
~
    \procname{insert}(\REF\LEXIKON \dictptr,\ \BITSTR x,\ \VALUE\ v) \idgl
    \BEGIN
        \REF\lexkonf \findpos,\insertpos,\newpos;
        \NATN $d$; \BOOL\ \keyexists;
~
        \findpos \gets \dictptr\REF\bdot\gate;
        \label{r11}%
        \WHILE(\findpos \neq \nil) \DO \untab
        \BEGIN
            \IF (\findpos < \dictptr) \untab
            \THEN \findpos \gets \ccopy(\findpos\REF,\ \dictptr\REFD\freeloc);%
            \untab
            \keyexists \gets (\findpos\REF\bdot\rkey\REF = x);
            \IF (\keyexists) \THEN \BREAK; \untab \untab
            \IF %
             (\text{ one of } x \text{ and } \findpos\REF\bdot\rkey\REF %
                \text{ is a left substring of the other }) \untab
            \THEN d \gets %
                \min(\length(x),\length(\findpos\REF\bdot\rkey\REF)) + 1;
            \ELSE \keyword{let}\yspace d \text{ be the smallest } d, %
                \text{ such that }%
                x[d] \neq (\findpos\REF\bdot\rkey\REF)[d]; % \untab
            \UNDENT
        \begin{minipage}^{\noobeycr}[t]{.87\textwidth}
            /* \\
             Note that $\findpos < \dictptr$ iff \findpos
             has been created before last call of
             \procname{open_environment}, which caused
             the creation of \dictptr.
             If this is true, then \findpos shall be replaced
             by a copy of it which is located atop of \dictptr.
             This copy will be removed at next call of
             \procname{close_environment}. \\ %
            **/
        \end{minipage}
            \insertpos \gets \findpos; %
            \findpos \gets \findpos\REF\bdot\zweig;
            \IF (\findpos \neq \nil \land \findpos < \dictptr) \THEN %
            \untab \untab
            \BEGIN
                \newpos \gets \ccopy(\findpos\REF,\ \dictptr\REFD\freeloc);
                \insertpos\REF\bdot\zweig \gets \newpos;
                \findpos \gets \newpos;
            \END;
            \WHILE (\findpos \neq \nil \land d \gt \findpos\REF\bdot\dist) %
            \DO \untab
            \BEGIN
                \insertpos \gets \findpos;
                \findpos \gets \findpos\REF\bdot\link;
                \IF (\findpos \neq \nil \land \findpos < \dictptr) \THEN%
                \untab\untab
                \BEGIN
                    \newpos \gets \ccopy(\findpos\REF,\ \dictptr\REFD\freeloc);
                    \insertpos\REF\bdot\link \gets \newpos;
                    \findpos \gets \newpos;
                \END;
            \END;
            \label{r12}% \Delta:\enspace %
            \IF (\findpos = \nil \lor d \neq \findpos\REF\bdot\dist) %
                \THEN \BREAK; \untab \untab
            \mbox{/* with \keyexists=\false, i.e. $x$ is not a valid key */}
~
            x \gets \enangle{x[d+j]}_{j<\length(x)-d}; %
            \mbox{/* $\idgl \textsf{shift_right}(x,d)$ */}
            \mbox{/* searching continues in branching off %
                  (after position $d$ of $\rkey\REF$) */}
        \END;
~
        \IF (\keyexists) \untab
        \THEN \findpos\REF\bdot\val \gets v
        \ELSE \untab
        \BEGIN
            \newpos \gets \create(\LEXKONF, \dictptr\REF\bdot\freeloc);
            \newpos\REFD\dist \gets d;
            \IF (\insertpos = \nil) %
            \THEN \dictptr\REF\bdot\gate \gets \newpos; \untab \untab
            \mbox{/* note that } \insertpos=\nil %
            \mbox{ iff } \dictptr\REF\bdot\gate = \nil \mbox{ */}
            \IF (\insertpos\REFD\branch = \findpos) \untab
            \THEN \insertpos\REFD\branch \gets \newpos
            \ELSE \insertpos\REFD\link \gets \newpos; \untab
            \newpos\REFD\branch \gets \nil;
            \newpos\REFD\link \gets \findpos;
            \newpos\REFD\rkey \gets %
            \ccopy(\enangle{x[d+j]}_{j<\length(x)-d},
\dictptr\REF\bdot\freeloc);
            \newpos\REFD\val \gets v;
        \END;
    \END. \mbox{/* \textit{insert} */}
\end{program}
\end{DATASTR}

There is a stack of items of type \LEXIKON with top \dictptr.
These items are connected by the downward link \fieldsel{pop}.
Next segment below \dictptr\REF\ is \dictptr\REFD\fieldsel{pop}\REF.
Each element (of type \LEXKONF or \BITSTR) we can access from
\dictptr\REFD\fieldsel{pop}\REF\ was created before \dictptr\REF.
By Definition of \create\ and \ccopy\ an element $p\REF$ is older
(created later) than another element $q\REF$ if and only if $p < q$.
Hence the address value of each element accessible from
\dictptr\REFD\fieldsel{pop}\REF\
(by iterated application of \bdot\zweig\REF, \bdot\link\REF)
is less than \dictptr.
%% The stack with top \dictptr
\par
Part \tit{39-73}\ of procedure \procname{insert}\
is a modification of \procname{search}\ from \ref{lexikon},
which additionally copies the traversed entries \findpos\REF\
and chain the copies parallel to the  originals by
\zweig\ and \link.
The same considerations as for \procname{search}\ yield that
the worst case time complexity of procedure \procname{insert}\
also is $\mathsf{O}(s)$.
%% \END{lifo03}
%    \vfill \small \bibliography{lifo}
\end{document}